\documentclass{article}
\usepackage[french]{babel}

\usepackage[T1]{fontenc}
\usepackage{graphicx,amsmath,amssymb,enumitem,tikz,footnote}
\usetikzlibrary{quotes,angles}

\usepackage[letterpaper,top=2.5cm,bottom=2.5cm,left=3cm,right=3cm,marginparwidth=1.75cm]{geometry}

\usepackage[colorlinks=true, allcolors=blue]{hyperref}

\title{Marches al\'eatoires dans un c\^one et fonctions discrètes harmoniques}

\author{Kilian Raschel et Pierre Tarrago}

\begin{document}
\maketitle

\begin{abstract}
Les marches aléatoires dans un cône présentent le double attrait de se trouver au c\oe ur de nombreux problèmes probabilistes et d'être liées à de multiples domaines mathématiques, comme la théorie spectrale, la combinatoire ou l'analyse complexe discrète. Dans cet article nous présentons quelques éléments-clé associés à ces processus : nous évoquerons leur définition, le lien avec le mouvement brownien dans un cône, ainsi que quelques sujets de recherche récents comme la construction de fonctions harmoniques discrètes.
\end{abstract}

\section{Introduction}
\label{sec:introduction}

\subsection{Marches aléatoires et comportement asymptotique}
Dans cet article, nous nous intéresserons au comportement de \textit{marches aléatoires} dans des domaines coniques et commencerons par présenter les motivations et contours de ce sujet. Rappelons qu'une marche aléatoire $S=(S_n)_{n\geq1}$ est un processus aléatoire à temps discret évoluant dans un espace d'états $E$ (que nous considérerons ici toujours dénombrable), tel que la position de $S$ au temps $n+1$ dépend du passé uniquement à travers sa position au temps $n$. Formellement, nous avons donc 
\begin{equation}
\label{eq:markov}
   \mathbb{P}(S_{n+1}=y\vert S_m=x_m,1\leq m\leq n)=\mathbb{P}(S_{n+1}=y\vert S_n=x_n):=k_{n}(x_n,y),
\end{equation}
pour tous $y$ et $x_1,\ldots,x_n$ appartenant à $E$, où $\mathbb{P}(A\vert B)$ est la probabilité  d'un évènement $A$ sachant que l'évènement $B$ est réalisé, appelée probabilité de $A$ conditionnellement à $B$. Nous ferons l'hypothèse que la loi de transition $k:=k_n$ est indépendante du temps (on dit alors que la marche aléatoire est \textit{homogène en temps}). Une des questions principales est la compréhension de la marche aléatoire en temps long : sachant que la marche aléatoire part d'un point $x_1$ au temps $n=1$, on s'intéresse à
\begin{equation*}
   \mathbb{P}(S_n=y\vert S_1=x_1),
\end{equation*}
la probabilité qu'elle atteigne un point $y$ en un temps $n$ long.

Deux cas classiques se présentent alors. Si $E$ est fini, la situation est complètement comprise au premier ordre : sous certaines hypothèses génériques d'\textit{irréductibilité} et d'\textit{apériodicité}, ce qui signifie simplement que la marche aléatoire peut aller d'un point à l'autre de l'espace d'états sans évolution périodique, il existe une mesure de probabilité $m$ sur $E$ telle que quand $n\rightarrow \infty$,
\begin{equation*}
   \mathbb{P}(S_n=y\vert S_1=x_1)\rightarrow m(y)
\end{equation*}
pour tous $x_1,y\in E$. Observons en particulier qu'en temps long, la marche aléatoire oublie son point de départ.

Le deuxième cas classique est celui d'un réseau $E=\mathbb{Z}^d$, avec $d\geq 1$, et d'un \textit{noyau} $k$ invariant par translation: $k(x,y)=k(0,y-x):=f(y-x)$. Dans ce cas, qui est la situation principale analysée dans cet article, pour $x_0\in \mathbb{Z}^d$, on notera $(x_0+S_n)_{n\geq 0}$ la marche aléatoire décalée de $x_0$, avec $S_0=0$. Si la marche est vraiment $d$-dimensionelle (ce qui signifie que le support de $f$ génère un sous-module de $\mathbb{Z}^d$ de rang $d$), il existe alors une constante $c>0$ telle que
\begin{equation*}
   \mathbb{P}(S_n=y)\leq \frac{c}{n^{d/2}}.
\end{equation*}
Avec des hypothèses supplémentaires, il est possible d'avoir un développement asymptotique précis de la probabilité précédente, ce que l'on appelle un \textit{théorème de la limite locale}. Par exemple, si $f$ admet un moment d'ordre deux et une dérive nulle, c'est-à-dire si $\sum_{x\in \mathbb{Z}^d}f(x)\vert x\vert^2<\infty$ et $\sum_{x\in \mathbb{Z}^d}f(x)x=0$, alors 
\begin{equation}
\label{eq:TLL_non-contraint}
   \mathbb{P}(S_n=y)\sim \frac{1}{(2\pi n)^{d/2}\Gamma},
\end{equation}
où $\Gamma$ est lié aux moments d'ordre deux de $f$. Une fois encore, la marche aléatoire oublie son point de départ en première approximation. 

L'\textit{invariance par translation} du noyau $k$ est fondamentale pour l'obtention de ces comportements locaux. Il existe d'ailleurs des généralisations de ces résultats à des marches sur des groupes plus généraux que $\mathbb{Z}^d$ ou sur des espaces homogènes, où une telle invariance par translation peut de nouveau être introduite. Au contraire, il est souvent très difficile d'obtenir un résultat asymptotique sur le comportement d'une marche aléatoire sur un espace qui soit à la fois infini et non invariant par translation; c'est le programme de ce qui suit dans le cas de cônes.

\begin{figure}
\begin{center}
\includegraphics[width=3.5cm]{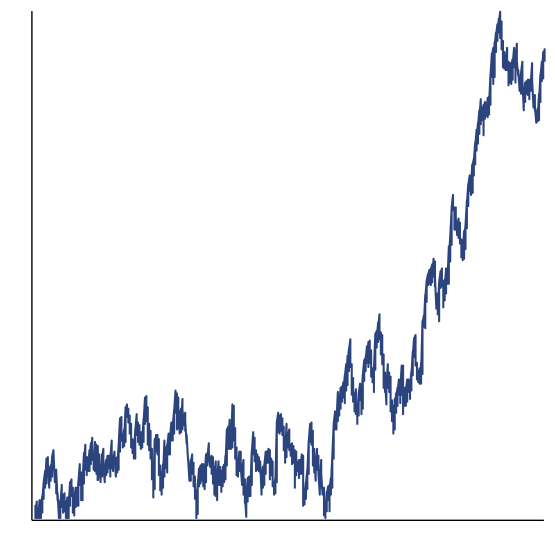}\qquad\qquad
\includegraphics[width=3.5cm]{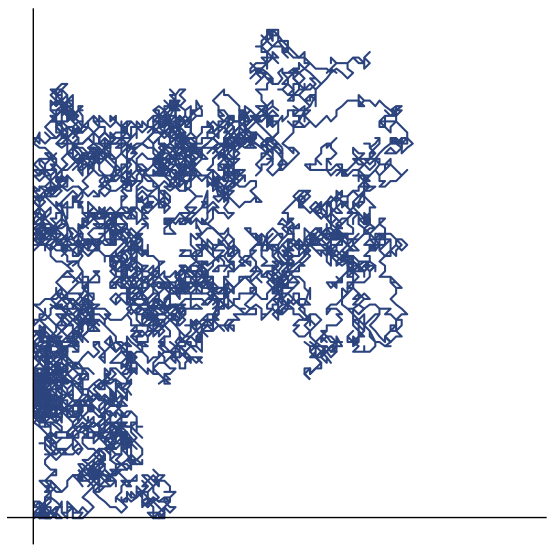}\qquad\qquad
\includegraphics[width=3.5cm]{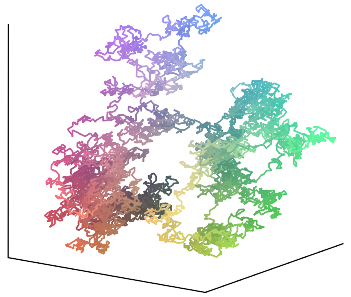}
\end{center}
\caption{Marches dans des cônes en petite dimension}
\label{fig:ex}
\end{figure}

\subsection{Marches aléatoires confinées}
\label{subsec:marche_conditionnee}

Un moyen simple de considérer des espaces d'états infinis et non invariants par translation est de restreindre une marche aléatoire sur $\mathbb{Z}^d$ à un sous-ensemble $A$ infini du réseau. Sans vouloir faire écho à l'actualité sanitaire, il n'est pas simple d'imposer un confinement à notre processus: en effet, dans de nombreux cas $\mathbb{P}\bigl(\bigcap_{n\geq 0}\{x_0+S_n\in A\}\bigr)=0$. Comment conditionner par un évènement de mesure nulle?

Une façon élémentaire de contourner ce problème utilise le premier temps de sortie de~$A$
\begin{equation}
\label{eq:temps_de_sortie}
   \tau^A_{x_0}=\inf\{n\geq 0: x_0+S_n\not\in A\}.
\end{equation}
Nous pouvons introduire la marche $S$ conditionnée à rester dans le domaine $A$ jusqu'à un certain temps $N\geq 0$, c'est-à-dire $\tau^A_{x_0}> N$. Cela donne bien une marche aléatoire jusqu'au temps $N$ (c'est-à-dire que \eqref{eq:markov} est encore satisfaite par le processus après conditionnement), dont les probabilités de transition sont données, pour $x,y\in A$ et $n<N$, par 
\begin{equation}
\label{eq:conditionnement_bef}
   \mathbb{P}(x_0+S_{n+1}=y\vert x_0+S_n=x, \tau^A_{x_0}> N)=
   k(x,y)\frac{\mathbb{P}(\tau_{y}^A> N-n-1)}{\mathbb{P}(\tau_{x}^A>N-n)}.
\end{equation}
Deux inconvénients sont intrinsèques à la construction précédente : tout d'abord ce nouveau processus n'est défini que jusqu'au temps $N$; de plus, la nouvelle marche n'est homogène ni en temps ni en espace, et cela même si la marche initiale possédait cette double homogénéité.

Un moyen de retrouver une homogénéité en temps et un horizon de temps infini est de considérer la limite de \eqref{eq:conditionnement_bef} quand $N$ tend vers l'infini. On obtient alors un nouveau processus $S^A=(S^A_n)_{n\geq0}$ dont les transitions valent, pour $x,y\in A$,
\begin{equation}
\label{eq:conditionnement}
   \mathbb{P}(x_0+S^A_{n+1}=y\vert x_0+S^A_{n}=x)=k(x,y)\lim_{N\rightarrow \infty}\frac{\mathbb{P}(\tau_{y}^A> N-1)}{\mathbb{P}(\tau_{x}^A>N)}.
\end{equation}
Sous réserve d'existence de la limite ci-dessus (ce n'est en rien acquis!), le nouveau processus $S^A$ définit une marche aléatoire  homogène en temps et à valeurs dans $A$ (par construction). Son noyau est pourvu d'une structure bien particulière, qui fait un lien avec l'analyse harmonique. Fixons $x_0\in A$ et posons
\begin{equation}
\label{eq:notation_h}
   h(y)=\lim_{N\rightarrow\infty}\frac{\mathbb{P}(\tau_{y}^A> N-1)}{\mathbb{P}(\tau_{x_0}^A>N)}.
\end{equation}
Nous noterons $\lambda=h(x_0)$. On calcule alors la limite qui apparaît dans \eqref{eq:conditionnement}:
\begin{equation*}
   \lim_{N\rightarrow \infty}\frac{\mathbb{P}(\tau_{y}^A> N-1)}{\mathbb{P}(\tau_{x}^A>N)}=\lambda\frac{h(y)}{h(x)}.
\end{equation*}
La relation $\mathbb{P}(\tau_{x}^A> N)=\sum_{y\in A}k(x,y)\mathbb{P}(\tau^A_y>N-1)$, obtenue par la propriété de Markov, conduit à
\begin{equation}
\label{eq:harmonicite}
   h(x)=\lambda\sum_{y\in A}k(x,y)h(y).
\end{equation}
En d'autres termes, la fonction $h$ dans \eqref{eq:notation_h} est \textit{$\lambda$-harmonique discrète} pour le noyau $k$ et dans le domaine $A$.

Notre notation \eqref{eq:notation_h} permet de reformuler le noyau \eqref{eq:conditionnement} de la marche conditionnée $S^A$ comme
\begin{equation}
\label{eq:Doob_conditioning}
   k^A(x,y)=\lambda\frac{h(y)}{h(x)}k(x,y),
\end{equation}
pour $x,y\in A$. La transformation \eqref{eq:Doob_conditioning} d'un noyau $k$ par une fonction $h$ est appelée \textit{conditionnement de Doob}, voir \cite{Do-59}. Cette construction fournit un nouveau noyau de marche aléatoire seulement si $h$ est positive et $\lambda$-harmonique au sens de \eqref{eq:harmonicite}. Cette méthode, plus conceptuelle que l'approche intuitive présentée en \eqref{eq:conditionnement_bef}, résout le problème du conditionnement par un évènement de mesure nulle. Les fonctions harmoniques jouent donc un rôle central dans le conditionnement de marches aléatoires de $\mathbb{Z}^d$ à rester dans des sous-domaines fixés.

De nombreuses questions se posent cependant au sujet de la construction précédente:
\begin{enumerate}[label=(Q\arabic{*}),ref=(Q\arabic{*})]
   \item\label{Q:question_1}Dans quels cas a-t-on bien convergence du quotient des \textit{probabilités de survie} dans  \eqref{eq:conditionnement} ?
   \item\label{Q:question_2}Quand \eqref{eq:conditionnement} est bien défini, il existe donc une fonction $\lambda$-harmonique $h$ introduite en \eqref{eq:notation_h}, directement liée au comportement probabiliste de la marche aléatoire en temps long. Peut-on construire et interpréter d'autres fonctions $\lambda$-harmoniques, qui conduiraient à d'autres conditionnements de Doob ?
   \item\label{Q:question_3}Le conditionnement \eqref{eq:conditionnement} (ou \eqref{eq:Doob_conditioning}) de la marche aléatoire est entièrement décrit par la fonction harmonique $h$. Existe-t-il, pour certaines marches aléatoires au moins, un moyen de calculer explicitement ou numériquement $h$ ?
\end{enumerate}
Ces trois questions, qui guideront notre présentation, restent largement ouvertes à l'heure actuelle. Le présent article va résumer les résultats récents obtenus dans le cas où le domaine de confinement $A$ est l'intersection d'un \textit{cône} $\mathcal{C}$ de $\mathbb{R}^d$ avec le réseau $\mathbb{Z}^d$. Tout en étant infinis et non invariants par translation, ces domaines sont omniprésents en probabilités pures et appliquées : ils sont par exemple régulièrement utilisés en théorie des files d'attente \cite{CoBo-83,FaIaMa-17}, en combinatoire \cite{BMMi-10} ou encore en théorie des représentations des groupes de Lie \cite{Bi-91}. Ils présentent de plus l'avantage d'être invariants par dilatation et d'être très bien compris d'un point de vue analytique: on sait par exemple décrire explicitement le noyau de la chaleur avec conditions de Dirichlet dans un cône $\mathcal{C}$ assez régulier \cite{Va-99}, ce qui s'avère essentiel pour l'étude des marches aléatoires conditionnées (voir Section \ref{sec:derive_nulle}).

Formellement, nous allons fixer un sous-domaine connexe et ouvert $\Sigma$ de la sphère unité $\mathbb{S}^{d-1}\subset \mathbb{R}^d$ et considérerons le cône
\begin{equation*}
   \mathcal{C}=\mathbb{Z}^d\cap (\mathbb{R}_+^*\cdot \Sigma).
\end{equation*}
Le domaine $A$ considéré est donc l'intersection de $\mathbb{Z}^d$ avec le cône engendré par l'ensemble des demi-droites issues de $0$ et passant par $\Sigma$.

\section{Marches homogènes à dérive nulle et unicité de la fonction harmonique}
\label{sec:derive_nulle}

Le cas le plus simple est celui d'une marche aléatoire homogène et à dérive nulle. En effet, pour notre cas conditionné on peut tirer profit de la très riche littérature sur les marches aléatoires à dérive nulle dans $\mathbb Z^d$ (sans contrainte de cônes). Le cas d'une dérive non nulle sera abordé en Section~\ref{sec:beyond}.

Pour simplifier la description, nous supposerons que les accroissements de la marche ont une matrice de covariance égale à l'identité, c'est-à-dire, pour une marche sans dérive,
\begin{equation*}
   \bigl(\mathbb{E}[S_{1}(i)S_1(j)]\bigr)_{1\leq i,j\leq d}
\end{equation*}
vaut la matrice identité, où $S_1(i)$ est la $i$-ème coordonnée de $S_1$ avec, rappelons-le, $S_0=0$. Par un changement linéaire approprié, les résultats de cette section restent toutefois valides sans cette restriction.

\subsection{Probabilité de survie}

Les nombreux résultats décrivant le comportement de marches dans des cônes culminent avec l'article \cite{DeWa-15} de Denisov et Wachtel: il en ressort notamment que sous des hypothèses très faibles, il existe un paramètre $p>0$ intrinsèque au cône $\mathcal C$ et une fonction $V$ positive et $1$-harmonique (au sens de \eqref{eq:harmonicite}) tels que le temps de sortie $\tau_{x}^{\mathcal C}$ (voir notre notation \eqref{eq:temps_de_sortie}) vérifie quand $N\rightarrow \infty$
\begin{equation}
\label{eq:tempslong_discret}
   \mathbb{P}(\tau_{x}^{\mathcal C}>N)\sim \frac{V(x)}{N^{p/2}}.
\end{equation}
Cette estimée assure en particulier la bonne définition de la marche conditionnée construite en \eqref{eq:conditionnement}, ce qui répond à la question \ref{Q:question_1} : pour tout choix de $x_0\in \mathbb{Z}^d\cap \mathcal{C}$, la fonction $h$ définie en \eqref{eq:notation_h} est donc égale à $\frac{V}{V(x_0)}$ et $\lambda=h(x_0)=1$. Dans la suite, nous dirons  simplement qu'une fonction est harmonique quand elle est $1$-harmonique.

\subsection{Du brownien aux marches}

Nous souhaitons à présent mentionner un outil-clé utilisé par les auteurs de \cite{DeWa-15} pour prouver \eqref{eq:tempslong_discret}, et plus généralement pour aboutir à une description fine de la marche aléatoire conditionnée en temps long. Il s'agit d'un \textit{couplage} entre marche aléatoire et \textit{mouvement brownien}.

Rappelons qu'un mouvement brownien $(B_t)_{t\geq 0}$ $d$-dimensionnel est un chemin aléatoire continu dans $\mathbb{R}^d$ issu de $0$; il représente l'objet limite universel des marches aléatoires sans dérive 
et non conditionnées: après \textit{renormalisation} conjointe des trajectoires de $S_n$ par $\sqrt{N}$ entre les instants $1$ et $N$ et du temps par $N$, la trajectoire discrète ressemble à un mouvement brownien $(B_t)_{0\leq t\leq 1}$ quand $N$ devient grand. Cela implique en particulier qu'entre les temps $1$ et $N$, l'ordre de fluctuation de ce type de marches aléatoires est de l'ordre $\sqrt{N}$.

Un phénomène plus précis a lieu: sans même renormaliser la marche aléatoire, il existe un mouvement brownien qui reste en tout temps $1\leq n\leq N$ à distance négligeable de la marche aléatoire par rapport à l'ordre des fluctuations. En particulier, si ce mouvement brownien est loin des bords du cône, avec grande probabilité la marche aléatoire couplée s'en tiendra également éloignée. Ceci est très important pour notre étude : au lieu de considérer une marche discrète dans un domaine restreint, domaine de recherche pour lequel on a peu d'outils directs, il suffit de considérer un mouvement brownien pour lequel on a beaucoup plus d'informations, comme en témoigne par exemple l'un des articles fondateurs \cite{De-87} de ce domaine.

Le comportement de la diffusion brownienne dans un domaine $\mathcal{C}$ de $\mathbb{R}^d$ est plus simple à décrire que celui d'une marche aléatoire dans le même domaine : la (densité de la) probabilité $P_{t}(x,y)$ qu'un brownien issu de $x\in \mathcal{C}$ se trouve au voisinage de $y\in \mathcal{C}$ au temps $t$ sans avoir quitté $\mathcal{C}$ suit l'équation de la chaleur avec conditions de Dirichlet au bord $\partial \mathcal{C}$ de $\mathcal{C}$, c'est-à-dire 
\begin{equation*}
   \left\lbrace\begin{aligned}
\partial_{t}P_{t}(x,\cdot)-\Delta P_{t}(x,\cdot)&=0,\quad  \text{sur }\mathcal{C},\\
P_{t}(x,y)&=0 ,\quad y\in\partial \mathcal{C},\\
\lim_{t\rightarrow 0} P_{t}(x,\cdot)&=\delta_{x},
\end{aligned}\right.
\end{equation*}
où $\lim_{t\rightarrow 0} P_{t}(x,\cdot)$ représente la limite au sens faible de la suite de fonctions $(P_{t}(x,\cdot))_{t>0}$, $\delta_x$ est la fonction de Dirac en $x$, et $\Delta$ désigne le Laplacien classique de $\mathbb R^d$. De même, la probabilité de survie $P_t(x)$ (probabilité que le mouvement brownien issu de $x$ n'ait pas quitté $\mathcal{C}$ à l'instant $t$) obéit à l'équation de la chaleur avec conditions de Dirichlet au bord et condition initiale $P_0(x)=1$, $x\in\mathcal{C}$. 

L'équation de la chaleur dans un cône est classique du point de vue de l'analyse spectrale; le comportement des solutions de ces équations aux dérivées partielles est gouverné par une fonction $u:\mathcal{C}\rightarrow\mathbb{R}_+$ qui satisfait l'équation d'harmonicité $\Delta u=0$ sur $\mathcal{C}$ et $u_{\partial \mathcal{C}}=0$. Il se trouve que cette fonction, appelée la {\it réduite} du cône, est l'unique fonction harmonique positive sur $\mathcal{C}$ s'annulant sur le bord (à une constante multiplicative près). Pour donner une idée du rôle fondamental joué par cette fonction dans le cas discret, notons que la fonction harmonique $V$ donnée en \eqref{eq:tempslong_discret} admet le comportement asymptotique suivant, quand $\vert x\vert \rightarrow \infty$ et $d(x,\partial \mathcal{C})\geq \varepsilon \vert x\vert$:
\begin{equation*}
   V(x)\sim u(x)
\end{equation*}
quel que soit $\varepsilon >0$, voir \cite{DeWa-15}. La fonction $V$ étant un exemple de fonction harmonique $h$ introduite en Section~\ref{sec:introduction}, cette asymptotique répond partiellement à \ref{Q:question_3} loin des bords de $\mathcal{C}$. En particulier, quand $x$ et $y$ sont grands et se tiennent éloignés du bord du cône, la marche conditionnée $S^A$ dans le domaine $A=\mathcal{C}\cap \mathbb{Z}^d$ grâce à la fonction harmonique $h=\frac{V}{V(x_0)}$ a un noyau \eqref{eq:Doob_conditioning} asymptotiquement égal à 
\begin{equation*}
   k^A(x,y)\sim \frac{u(y)}{u(x)}k(x,y).
\end{equation*}

\subsection{Théorie de Doob et unicité de la fonction harmonique}
\label{subsec:unicite}

Les résultats de Denisov et Wachtel apportent une réponse positive à \ref{Q:question_1} en s'appuyant sur la réduite $u$ du cône $\mathcal{C}$, qui est l'unique fonction harmonique positive pour l'opérateur Laplacien dans $\mathcal{C}$ satisfaisant aux conditions de Dirichlet. Il est donc tentant de supposer que de façon analogue, il existe une unique fonction harmonique discrète positive solution de \eqref{eq:harmonicite}, à des constantes multiplicatives près. Cependant, comme nous avons pu écrire plus haut, l'une des difficultés posée par les marches conditionnées consiste en l'absence d'outils analytiques directs (voir toutefois les Sections~\ref{subsec:Harnack} et \ref{sec:analytique} pour des approches analytiques possibles).

Il existe tout de même un moyen général de décrire l'ensemble des fonctions harmoniques positives d'une marche aléatoire, que nous allons maintenant aborder. Nous avons vu dans la Section~\ref{subsec:marche_conditionnee} que conditionner une marche sur $\mathbb{Z}^d$ à rester dans un domaine jusqu'à un temps long faisait apparaître, à la limite, une fonction harmonique particulière, voir \eqref{eq:harmonicite} : celle-ci se caractérise comme la limite d'un quotient de probabilités de survie, à condition que ce dernier converge. On souhaiterait modifier cette définition pour se débarrasser de la condition de convergence et ainsi obtenir toutes les fonctions harmoniques possibles. C'est précisément l'objectif de la théorie de la \textit{frontière de Martin}.

\subsubsection*{Frontière de Martin}
\label{subsec:Martin}

Cette théorie, dont nous ne pouvons donner qu'un court aperçu ici, relie de manière très précise l'ensemble des fonctions harmoniques positives au comportement asymptotique de la marche aléatoire sous-jacente. Les références \cite{Do-59,KuMa-98,IR-08} contiennent une présentation détaillée de cette notion.

Élargissons d'abord notre problème en nous intéressant à l'équation 
\begin{equation}
\label{eq:super_harmonic}
\left\lbrace\begin{array}{ll}
\displaystyle g(x)-\sum_{y\in \mathcal{C}}k(x,y)g(y)=f(x),& x\in \mathcal{C},\\
\displaystyle g(x)=0,& x\in \partial \mathcal{C},
\end{array}\right.
\end{equation}
où $f:\mathcal{C}\rightarrow \mathbb{R}_+$ est une fonction donnée. L'équation d'harmonicité \eqref{eq:harmonicite} est un cas particulier de \eqref{eq:super_harmonic}, posant $f=0$. Cela revient naïvement à inverser l'opérateur $\text{Id}-k_{\,\vert \mathcal{C}}$ (qui, rappelons-le, n'est pas injectif puisque la fonction $V$ introduite en \eqref{eq:tempslong_discret} est dans son noyau); on voudrait donc poser 
\begin{align}
g(x)=&\sum_{y\in \mathcal{C}}\left(\sum_{n\geq 0}k_{\,\vert \mathcal C}^n(x,y)\right)f(y)\nonumber\\
=&\sum_{y\in \mathcal{C}}\left(\sum_{n\geq 0}\mathbb P(x+S_n=y,\tau_x^{\mathcal C}>n)\right)f(y):=Gf(x).\label{eq:notation_Gf}
\end{align}
Cette solution est bien définie pour toute fonction $f$ à support fini si et seulement si la \textit{fonction de Green}
\begin{equation}
\label{eq:def_Green}
   G(x,y):=\sum_{n\geq 0}\mathbb P(x+S_n=y,\tau_x^{\mathcal C}>n)
\end{equation}
est finie pour tous $x,y\in \mathcal{C}$. Cette condition de finitude est réalisée dès lors la marche aléatoire est \textit{transiente}, ce qui signifie qu'elle ne passe qu'un nombre fini de fois par chaque état. 

L'hypothèse de transience est bien satisfaite dans notre situation : du fait de la dérive nulle, la marche aléatoire non conditionnée sort forcément du cône $\mathcal{C}$ en temps fini.

On vérifie facilement que la fonction de Green $G(\cdot,y)$ définie en \eqref{eq:def_Green} est solution de l'équation \eqref{eq:super_harmonic} pour $f=\delta_y$. En d'autres termes, elle est harmonique partout sauf en $y$ ! L'observation-clé pour construire une fonction harmonique \textit{en tout point} et positive réside dans la remarque suivante : prenons une suite de fonctions $(f_n)_{n\geq0}$ de $\mathcal{C}$ dans $\mathbb{R}_+$ telle que
\begin{itemize}
   \item[$\bullet$]quand $n\to\infty$, $f_n(x)\to 0$ pour tout $x\in\mathcal{C}$;
   \item[$\bullet$]$(Gf_n)_{n\geq 0}$ converge simplement sur $\mathcal{C}$ vers une fonction non nulle (les fonctions $Gf_n$ étant définies par \eqref{eq:notation_Gf}).
\end{itemize}
Sous ces conditions, $\lim_{n\rightarrow\infty} Gf_n$ sera une fonction harmonique positive.

L'exemple fondamental d'une telle construction est le suivant : fixons $x_0\in \mathcal{C}$ et prenons
\begin{equation}
\label{eq;def_f_n_dirac}
   f_n=\frac{1}{G(x_0,y_n)}\delta_{y_n},
\end{equation}
pour une suite de points $(y_n)_{n\geq0}$ allant vers l'infini dans $\mathcal{C}$ telle que la suite $(Gf_n)_{n\geq 0}$ converge (alors nécessairement vers une fonction non nulle, puisque prenant la valeur $1$ en $x_0$). La limite définit une fonction harmonique, qui dépend \textit{a priori} de la façon avec laquelle la suite $(y_n)_{n\geq0}$ va à l'infini, mais ne dépend de $x_0$ qu'à une constante multiplicative près. L'ensemble des fonctions harmoniques ainsi obtenues est appelé la frontière de Martin de la marche.

Le \textit{théorème de représentation de Martin} affirme alors que \textit{toute} fonction harmonique égale à $1$ en $x_0$ peut être obtenue comme combinaison convexe  d'éléments de la frontière de Martin (le théorème donne en plus une structure particulière, voir \cite{Do-59}, qui va au-delà du cadre de cette introduction). Ce théorème donne donc une description de l'ensemble des fonctions harmoniques positives, à condition de pouvoir décrire la frontière de Martin.

Présentons une application de cette théorie au problème de l'unicité de la fonction harmonique positive. Imaginons que pour $y_n$ tendant vers l'infini, la fonction de Green \eqref{eq:def_Green} admette l'asymptotique
\begin{equation}
\label{eq:decouplage}
   G(x,y_n)\sim \frac{1}{\vert x-y_n\vert^\alpha}V(x)F(y_n),
\end{equation}
où $V$ est la fonction harmonique donnée en \eqref{eq:tempslong_discret}, $F$ est une fonction indépendante de $x$ et $\alpha$ est un paramètre dépendant éventuellement de $y_n$. Pour la suite $f_n$ introduite en \eqref{eq;def_f_n_dirac} et en omettant tout problème de convergence quand $n$ tend vers l'infini, on aura
\begin{equation*}
   Gf_n(x)\sim \frac{V(x)F(y_n)}{V(x_0)F(y_n)}\frac{\vert x_0-y_n\vert^\alpha}{\vert x-y_n\vert^\alpha}\sim\frac{V(x)}{V(x_0)}.
\end{equation*}
On voit donc que pour toute suite $(y_n)_{n\geq 0}$ tendant vers l'infini dans le cône, $Gf_n$ converge vers la fonction harmonique positive $V$, à une constante multiplicative près. Cela prouve l'unicité de la fonction harmonique positive.

L'un des premiers travaux prouvant de tels résultats fut \cite{Bi-91}, en lien avec des algèbres de Lie. Dans le cas d'un cône $\mathcal{C}$ convexe ou assez régulier et sous certaines conditions de moments des accroissements de la marche aléatoire, une preuve du découplage \eqref{eq:decouplage} est donnée dans \cite{DuRaTaWa} et expliquée dans le prochain paragraphe.

\subsubsection*{Asymptotique du noyau de Green}
Pour comprendre le découplage asymptotique de $G(x,y)$ évoqué en \eqref{eq:decouplage}, il est instructif de décomposer la fonction de Green \eqref{eq:def_Green} comme suit: 
\begin{align}
   G(x,y)&=\sum_{n=o(\vert x-y\vert)^2}\mathbb P(x+S_n=y,\tau_x^{\mathcal C}>n) +\sum_{n\succeq \vert x-y\vert^2}\mathbb P(x+S_n=y,\tau_x^{\mathcal C}>n)\nonumber\\
   &:=G_{\text{negl}}(x,y)+G_{\text{fluc}}(x,y),\label{eq:simplification_green}
\end{align}
où le symbole ``$\succeq$'' désigne donc les indices $n$ absents de la première somme, de manière à n'en oublier aucun. En effet, nous savons qu'une marche aléatoire non conditionnée fluctue avec une amplitude $\sqrt{n}$ au temps $n$; cela rend naturelle la décomposition ci-dessus quand $\vert x-y\vert$ devient grand, et on s'attend à ce que $G_{\text{negl}}$ soit négligeable devant $G_{\text{fluc}}$. 

Intéressons-nous au second terme dans \eqref{eq:simplification_green}. Quand $n$ est nettement plus grand que $\vert x-y\vert$, la position relative de $x$ et $y$ importe peu; seules comptent les positions de $x$ et $y$ par rapport au bord. Le terme 
\begin{equation}
\label{eq:TLL}
   \mathbb P(x+S_n=y,\tau_x^{\mathcal C}>n)
\end{equation}
vaut donc approximativement 
\begin{equation*}
   Cn^{-d/2}\mathbb{P}(\tau_x^{\mathcal C}\geq n/2)\mathbb{P}(\overline{\tau}_{y}^{\mathcal C}\geq n/2),
\end{equation*}
où $C>0$ est une constante et $\overline{\tau}_{y}^{\mathcal C}$ représente le temps de sortie pour la marche $(-S_n)_{n\geq1}$ : en effet, du point de vue de $y$ la marche aléatoire va à l'envers ! La contribution $n^{-d/2}$ est la contribution typique venant du théorème central de la limite locale en dimension $d\geq 1$, voir \eqref{eq:TLL_non-contraint}. Utilisant \eqref{eq:tempslong_discret}, on déduit que lorsque $y\to\infty$,
\begin{equation*}
   G_{\text{fluc}}(x,y) \sim V(x)\sum_{n\succeq \vert x-y\vert^2}n^{-d/2-p/2}\mathbb P(\overline{\tau}_{y}^{\mathcal C}\geq n/2).
\end{equation*}

Quand $y$ tend vers l'infini dans le cône, deux situations peuvent se présenter : soit la distance de $y$ au bord reste macroscopique devant $\vert y\vert$, soit $y$ s'approche du bord de $\mathcal{C}$ en allant vers l'infini. Dans le premier cas, la marche renversée issue de $y$ ne voit pas le bord, et on peut ignorer le terme $\mathbb P(\overline{\tau}_{y}^{\mathcal C}\geq n/2)$. Dans le deuxième cas, la marche aléatoire partant de $y$ ne voit pas l'origine du cône mais un cône $\mathcal{C}_y$ tangent à $\mathcal{C}$ en un point proche de $y$. En raffinant les résultats de \cite{DeWa-15}, il est possible de montrer que $\overline{\tau}_{y}^{\mathcal C}$ se comporte de manière semblable à \eqref{eq:tempslong_discret} pour le cône tangent $\mathcal{C}_y$, ce qui correspond à l'asymptotique 
\begin{equation*}
   \mathbb{P}(\overline{\tau}_y^{\mathcal{C}}>n)\sim\mathbb{P}(\overline{\tau}_y^{\mathcal{C}_y}>n)\sim\frac{F(y)}{n^{q/2}},
\end{equation*}
où $F(y)$ ne tend pas vers $0$ et $q$ est le paramètre $p$ du cône tangent $\mathcal{C}_y$.
 En sommant les contributions de $G_{\text{fluc}}$ sur $n\succeq \vert x-y\vert^2$, on a alors le comportement asymptotique 
\begin{equation}
\label{eq:fluct_asymp}
   G_{\text{fluc}}(x,y)\sim V(x)F(y)\vert x-y\vert^{-d-p-q+2}.
\end{equation}

Le problème restant est de savoir quand, dans \eqref{eq:simplification_green}, $G_{\text{negl}}(x,y)$ est effectivement négligeable par rapport à $G_{\text{fluc}}(x,y)$. Posons
\begin{equation*}
   q_0=\sup_{\mathcal{K}\text{ cône tangent à $\mathcal{C}$}} q_{\mathcal{K}}.
\end{equation*}
Le raisonnement précédent via \eqref{eq:fluct_asymp} permet d'interpréter $d+p+q_0-2$ comme la décroissance maximale de $G_{\text{fluc}}(x,y)$ quand $y$ tend vers l'infini. Si $\mathbb{E}[\vert  S_1\vert^{d+p+q_0-2}]=\infty$, on peut construire un exemple (semblable aux contre-exemples pour l'inégalité de Markov quand la condition de moment n'est plus vérifiée) pour lequel 
\begin{equation*}
   G_{\text{negl}}(x,y)\geq\mathbb{P}(x+S_1=y ,\tau_{x}^\mathcal{C}>1)\geq \vert x-y\vert^{-d-p-q_0+2}
\end{equation*}
quand $y$ tend vers l'infini, et donc pour lequel, en vertu de \eqref{eq:fluct_asymp}, $G_{\text{negl}}(x,y)$ n'est pas négligeable devant $G_{\text{fluc}}(x,y)$. 

On peut montrer qu'une condition suffisante pour garantir la négligeabilité de $G_{\text{negl}}(x,y)$ consiste à imposer $\mathbb{E}[\vert  S_1\vert^{d+p+q_0-2+\varepsilon}]<\infty$ pour un $\varepsilon>0$ quelconque. On a alors, sous cette nouvelle condition de moment,
\begin{equation*}
   G(x,y) \sim V(x)F(y)\vert y-x\vert^{-d-p-q+2}
\end{equation*}
quand $\vert y\vert$ tend vers l'infini et, d'après le paragraphe sur la frontière de Martin, ce découplage entraîne l'unicité (à constante multiplicative près) de la fonction harmonique positive. 

\section{Au delà du cas homogène à dérive nulle}
\label{sec:beyond}

Nous allons présenter deux résultats importants, décrivant l'ensemble des fonctions harmoniques positives dans un cône dans le cas où la marche aléatoire a une dérive non nulle ou n'est pas homogène. Le cas d'une dérive non nulle, qui est traité en premier en suivant \cite{IR-08}, conduit à une caractérisation topologique de la frontière de Martin à l'aide d'une courbe particulière. Dans le cas d'une marche aléatoire non homogène avec dérive nulle, qui est abordé dans un deuxième temps et suit \cite{BoMuSi-15}, il est possible de retrouver l'unicité de la fonction harmonique, comme dans le cas homogène.

\subsection{Dérive non nulle: heuristique et résultats}

Commençons par évoquer le fameux théorème de Ney et Spitzer \cite{NeSp-66}, traitant du cas classique de marches aléatoires dans $\mathbb Z^d$ homogènes, irréductibles et avec dérive, pour lequel les fonctions exponentielles jouent un rôle central. Pour plus de commodité nous nous placerons en dimension $d=2$; les techniques et résultats se généralisent \textit{verbatim} en dimension quelconque.

Posons $L(u,v)=\sum_{(i,j)\in\mathbb Z^2} f(i,j)\exp(ui+vj)$ pour $(u,v)\in\mathbb{R}^2$ où, rappelons-le (voir la Section~\ref{sec:introduction}), $f(i,j)$ représente la probabilité de saut $\mathbb{P}(S_1=(i,j)\vert S_0=(0,0))$. Un calcul simple montre qu'une fonction exponentielle $h(i,j)=\exp(ui+vj)$ est harmonique si et seulement si son paramètre $(u,v)\in\mathbb R^2$ satisfait l'équation 
\begin{equation}
\label{eq:courbe_Ney_Spitzer}
  L(u,v)= 1.
\end{equation}
La fonction $L$ est convexe et, en $(u,v)=(0,0)$, $L(0,0)=1$, $\nabla L(0,0)=\mathbb{E}[S_1]\in \mathbb{R}^2$ et $\text{Hess}[L](0,0)$ est la matrice de covariance de $S_1$, donc définie positive.

\begin{figure}
\begin{center}
\includegraphics[width=4cm]{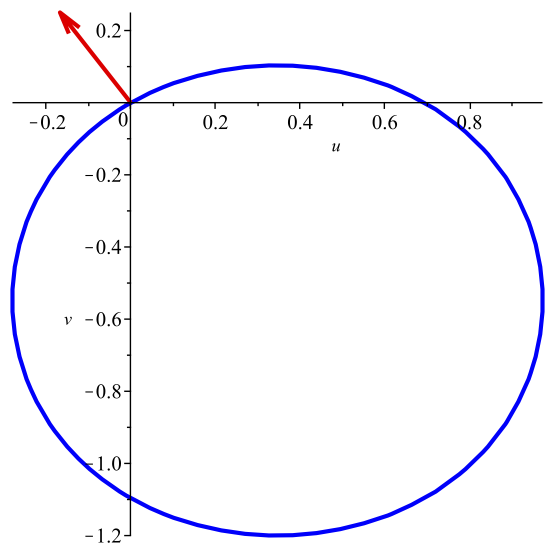}
\end{center}
\caption{Un exemple de courbe de niveau \eqref{eq:courbe_Ney_Spitzer}}
\label{fig:Hennequin}
\end{figure}

Une dérive nulle garantit pour une marche irréductible que seul le point $(u,v)=(0,0)$ est solution de \eqref{eq:courbe_Ney_Spitzer}, comme alors $(0,0)$ est le minimum global de $L$. \`A l'inverse, en présence d'une dérive, $L(0,0)$ n'est plus le minimum de $L$ et cette même équation \eqref{eq:courbe_Ney_Spitzer} est satisfaite par le bord d'un convexe fermé d'intérieur non vide. Elle définit donc un ensemble homéomorphe à $\mathbb S^{d-1}$, c'est-à-dire un cercle en dimension $2$. Nous représentons Figure~\ref{fig:Hennequin} la courbe de niveau \eqref{eq:courbe_Ney_Spitzer} pour la marche $L(u,v) = \frac{e^u}{6}+\frac{e^{-u}}{3}+\frac{3e^v}{8}+\frac{e^{-v}}{8}$. La dérive $(-\frac{1}{6},\frac{1}{4})$ est orthogonale à la courbe en l'origine.

Ney et Spitzer \cite{NeSp-66} obtiennent l'asymptotique des fonctions de Green
\begin{equation*}
   G(x,y)=\sum_{n\geq 0}\mathbb P(x+S_n=y)
\end{equation*}
et en déduisent que l'ensemble des $(u,v)$ caractérisés par \eqref{eq:courbe_Ney_Spitzer} représente exactement la frontière de Martin.

Avant de passer au cas d'un cône (plus restreint géographiquement mais plus riche mathématiquement), évoquons brièvement la propriété de \textit{minimalité} de certaines fonctions harmoniques. Une fonction harmonique positive $h$ est dite minimale si les seules fonctions harmoniques $g$ vérifiant $0\leq g\leq h$ sont en fait proportionnelles à $h$.

D'après la théorie de la frontière de Martin et le résultat de Ney et Spitzer, les exponentielles $h(i,j)=\exp(ui+vj)$ définissent les seules fonctions harmoniques minimales. On peut retrouver cette propriété par un calcul direct. Si en effet $h(i,j)$ est harmonique, alors $h(i+1,j)$ l'est également. Sous une hypothèse non contraignante d'irréductibilité, on trouve facilement $h(i,j)\geq c h(i+1,j)$, où par exemple $c$ est la probabilité de visiter $(i+1,j)$ en partant de $(i,j)$ si celle-ci est non nulle, et donc par minimalité $h(i+1,j)=a h(i,j)$. Comme un raisonnement symétrique s'applique à $h(i,j+1)$, on aboutit à la propriété classique de morphisme caractérisant les fonctions exponentielles.

Examinons maintenant la situation des marches dans des cônes, à travers l'exemple du demi-plan $\mathbb Z\times\mathbb N$. Les marches aléatoires réfléchies ou tuées dans des demi-espaces $\mathbb Z^{d-1}\times\mathbb N$ sont des cas particuliers de \textit{processus Markov additifs} (parfois appelés Markov modulés): l'une des coordonnées (ici la première, qui vit sur $\mathbb Z$) est la partie additive, invariante par translation; la seconde coordonnée est une chaîne de Markov sur $\mathbb N$ (\textit{a priori} inhomogène en espace) qui ``module'' le premier processus.

La structure des fonctions harmoniques minimales est claire: on peut appliquer l'argument précédent à $h(i+1,j)$ et déduire que toute fonction minimale harmonique prend la forme
\begin{equation}
\label{eq:u_a_trouver}
   h(i,j)=\exp(ui) h(0,j).
\end{equation}
Cela nous permet d'effectuer une réduction de dimension. Supposons un instant $u$ connu dans \eqref{eq:u_a_trouver}; alors la relation d'harmonicité discrète pour $h(i,j)$ se transforme en une relation analogue mais un-dimensionnelle pour $h(0,j)$. On bascule alors dans la théorie discrète du potentiel sur $\mathbb Z$ ou $\mathbb N$, que nous ne développerons pas ici, mais pour laquelle de nombreux outils et résultats sont disponibles, et grâce à laquelle $h(0,j)$ peut être caractérisée. 

Le point subtil est de trouver les valeurs possibles pour $u$ dans \eqref{eq:u_a_trouver}. Suivant une idée originale de Foley et McDonald, on peut calculer de deux façons différentes la limite de
\begin{equation}
\label{eq:quotient_Green}
   \frac{G\bigl((i+1,j),(k,\ell)\bigr)}{G\bigl((i,j),(k,\ell)\bigr)}
\end{equation}
lorsque $(k,\ell)$ tend vers l'infini dans la direction de la dérive. Choisir cette direction particulière n'est en fait pas une restriction: on peut toujours s'y ramener grâce à une transformation de Doob comme présentée en Section~\ref{sec:introduction}. D'une part, en utilisant les asymptotiques de chacune des fonctions de Green au numérateur et dénominateur, on aboutit à une limite dans \eqref{eq:quotient_Green} de la forme $h(i+1,j)/h(i,j)$, où $h$ est une certaine fonction harmonique, dont on peut montrer qu'elle est minimale. Par ailleurs, en utilisant une méthode de ``décomposition de Bernoulli'', qui étudie la dépendance de la fonction de Green en les premiers pas, toujours dans la direction asymptotique dictée par la dérive, on trouve que la limite de \eqref{eq:quotient_Green} est égale à $1$. Autrement dit, après conditionnement de Doob, on trouve systématiquement $u=0$ dans \eqref{eq:u_a_trouver}. 

Dans le cas des demi-espaces, les calculs sont menés dans \cite{IR-08}. Des idées similaires fonctionnent pour d'autres cônes, avec des détails plus complexes selon la nature des cônes.

\subsection{Marches inhomogènes et inégalités de Harnack}
\label{subsec:Harnack}

Les marches aléatoires $(x+S_n)_{n\geq0}$ considérées juqu'ici étaient homogènes en espace. Cela signifie (Section~\ref{sec:introduction}) que le noyau de la marche non conditionnée est de la forme $k(x,y)=f(y-x)$, où $f:\mathbb{Z}^d\rightarrow \mathbb{R}_+$ est une mesure de probabilité. Cette hypothèse est cruciale dans la Section~\ref{sec:derive_nulle}, puisqu'elle permet de réaliser un couplage entre marche aléatoire et mouvement brownien, expliquant ainsi les similitudes entre le comportement des fonctions harmoniques de la marche discrète et celles du cas diffusif. Cependant, une marche aléatoire générique est inhomogène en espace, et dans cette section nous présentons quelques idées qui s'appliquent dans ce cas non homogène. 

Répondre à \ref{Q:question_1} (asymptotique de la probabilité de survie comme en \eqref{eq:tempslong_discret}) dans le cas non homogène s'avère particulièrement complexe, en raison du comportement non régulier de la marche aléatoire. En revanche, il existe une approche permettant d'aborder la question \ref{Q:question_2} (description de l'ensemble des fonctions harmoniques). Cette méthode, décrite dans \cite{BoMuSi-15}, conduit au premier résultat général d'unicité de la fonction harmonique positive.

De nouveau, la méthode consiste à comparer la situation discrète engendrée par la marche aléatoire à la situation continue d'un processus de diffusion. Rappelons qu'une fonction harmonique pour le mouvement brownien est une fonction harmonique au sens (classique) analytique du terme, c'est-à-dire vérifiant dans le domaine
\begin{equation*}
   \Delta f=\sum_{i=1}^{d} \frac{\partial ^2f}{\partial x_i^2}=0.
\end{equation*}
Un comportement bien connu des fonctions harmoniques pour le Laplacien (valide en fait pour tout opérateur différentiel elliptique) est la propriété de \textit{moyennisation}, qui se traduit par l'\textit{inégalité de Harnack} : il existe une constante universelle $C>0$ (dépendant uniquement de la dimension) telle que si $f$ est une fonction positive harmonique sur une boule de centre $x$ et de rayon $2R$, alors pour tous $y,z\in B(x,R)$ on a 
\begin{equation}
\label{eq:Harnack}
   f(z)\leq Cf(y).
\end{equation}
Cette inégalité, qui peut s'obtenir directement à l'aide du noyau de Poisson dans le cas du Laplacien, se comprend intuitivement par le fait qu'une fonction harmonique varie peu quand elle s'éloigne du bord de son domaine de définition (précisément grâce à cette propriété de moyennisation).

Quand on impose de plus à la fonction harmonique de s'annuler sur le bord de son domaine, et si en outre ce bord est assez régulier (Lipschitz par exemple), on peut encore obtenir une inégalité de \textit{Harnack au bord}, qui prend la forme suivante. Il existe une constante universelle $C>0$ telle que pour toutes fonctions $f,g$ harmoniques positives sur l'intersection d'une boule $B(x,2R)$ et d'un domaine, si $f$ et $g$ s'annulent sur le bord de ce domaine, alors 
\begin{equation}
\label{eq:Harnack_boundary}
   \left\{\begin{array}{rcl}
   \frac{f(x)}{g(x)} &\leq& C\frac{f(e)}{g(e)},\smallskip\\
   f(x)&\leq& Cf(e),
   \end{array}\right.
\end{equation}
où $e$ est un point à distance macroscopique du bord (relativement à $R$). Par des méthodes assez générales, on peut montrer que les équations \eqref{eq:Harnack} et \eqref{eq:Harnack_boundary} assurent en même temps l'existence et l'unicité de la fonction harmonique positive s'annulant sur le bord, pour tout domaine au bord  assez régulier.

L'important travail mené dans \cite{BoMuSi-15} est d'obtenir ces inégalités de Harnack dans le cas de marches aléatoires non homogènes à dérive nulle dans l'orthant $\mathbb{N}^d$. Ces marches doivent tout de même satisfaire certaines conditions : tout d'abord elles doivent être elliptiques, ce qui signifie qu'il existe $\alpha>0$ tel que $k(x,x\pm e_i)\geq \alpha$ pour tout vecteur de la base canonique $e_i$. Cette condition est à rapprocher de la théorie des opérateurs différentiels elliptiques, pour lesquels les inégalités de Harnack existent. D'autre part, l'ensemble de pas doit être borné uniformément pour tout $x\in \mathbb{N}^d$. Cette dernière hypothèse est cruciale : dans le cas d'un ensemble de pas non bornés, les inégalités de Harnack peuvent ne plus être valables \cite{Law-93}. 

\`{A} l'aide de ces inégalités, les auteurs de \cite{BoMuSi-15} parviennent à prouver l'existence et l'unicité de la fonction harmonique positive (à constante multiplicative près) par des procédés très semblables à ceux utilisés dans le cas continu. En plus de constituer la première réponse à \ref{Q:question_2} dans un cadre général, ce résultat montre l'aspect central des inégalités de Harnack dans l'étude des chaînes de Markov discrètes. Cette approche a été poursuivie pour obtenir des estimées gaussiennes de marches aléatoires inhomogènes en espace dans des cônes, ou encore pour répondre à la question \ref{Q:question_2} dans le cas où le domaine sur lequel on restreint la marche aléatoire est un domaine Lipschitz.

\section{Sur la construction explicite des fonctions harmoniques dans un quadrant}
\label{sec:analytique}

Cette section est consacrée à la question~\ref{Q:question_3} dans le cas d'un quart de plan et d'une marche homogène.

Dans le cas classique continu, le calcul des fonctions harmoniques dans des cônes de $\mathbb R^d$ s'effectue facilement, grâce à la décomposition polaire. \`A titre d'exemple, en dimension $2$, toute fonction harmonique $f(x,y)=\widetilde{f}(\rho,t)$ dans le cône d'ouverture $\theta$ 
\begin{equation*}
   \{\rho e^{it} : \rho\geq 0 \text{ et } 0\leq t\leq \theta\}
\end{equation*}
s'écrira comme une somme pondérée de termes $\rho^{k\pi/\theta}\sin(kt\pi/\theta)$, avec $k\geq 1$. Dans le cas discret, l'absence de décomposition polaire sera source de complications considérables mais aussi de richesses.

La littérature existante donne des pistes sur des constructions possibles de fonctions discrètes harmoniques. Si des polynômes ou exponentielles harmoniques peuvent être obtenus de façon élémentaire, via une récurrence, c'est souvent l'existence d'une structure supplémentaire qui ouvre la voie à la construction de fonctions harmoniques. Par exemple, lorsque le cône est une chambre de Weyl, des outils comme le principe de réflexion, des déterminants de Vandermonde, ou un détour par la théorie des représentations et des réseaux de poids permettent d'exprimer des fonctions harmoniques, voir par exemple \cite{Bi-91}.

Nous souhaitons ici présenter une construction plus analytique, qui est l'objet du travail \cite{HoRaTa-20}. Notre approche, bien qu'intrinsèque à la dimension $2$ (elle fonctionne \textit{a priori} seulement dans le quart de plan -- ou de façon équivalente dans tout cône convexe planaire, après transformation linéaire), ne suppose aucune symétrie ou structure additionnelle sur les pas. L'idée se résume simplement : en utilisant des séries génératrices, les équations d'harmonicité discrète \eqref{eq:harmonicite} se reformulent en une \textit{équation fonctionnelle}, qui sera énoncée en \eqref{eq:main_functional_equation}. Celle-ci appartient à une classe d'équations qui, du fait de ses apparitions répétées dans des problèmes de files d'attente \cite{FaIaMa-17}, de processus stationnaires dans le quart de plan, de processus transients \cite{KuMa-98}, ou encore d'énumération des chemins en combinatoire \cite{BMMi-10}, s'est retrouvée au centre d'investigations poussées de la part de différentes communautés mathématiques. 

Notre construction permet de décrire l'ensemble des fonctions harmoniques réelles (non nécessairement positives), qui apparaissent notamment dans des développements asymptotiques de certaines quantités pertinentes, comme la probabilité de survie \eqref{eq:tempslong_discret}. Elle conduit donc à la description de la structure de l'ensemble des fonctions harmoniques.  
En outre, dans le cas où la frontière de Martin est un singleton (par exemple quand la dérive est nulle, comme prouvé dans \cite{DuRaTaWa}, voir notre Section~\ref{subsec:unicite}), notre construction permet également de calculer l'unique fonction harmonique positive, celle utilisée pour le conditionnement de Doob dans les Sections~\ref{sec:introduction} et \ref{sec:derive_nulle}.

\subsection{Une équation fonctionnelle...}
Afin d'énoncer notre équation principale, nous introduisons donc une marche aléatoire dans $\mathbb Z^2$ avec probabilité de saut $f(k,\ell)$ dans la direction $(k,\ell)$. Supposons que les sauts positifs sont petits, c'est-à-dire que $f(k,\ell)=0$ si $k>1$ ou $\ell>1$; en revanche, les sauts négatifs peuvent être choisis arbitrairement grands. Introduisons le noyau de la marche
\begin{equation*}
   K(x,y)=xy\bigl(\textstyle \sum_{k,\ell} f(k,\ell) x^{-k}y^{-\ell}-1\bigr).
\end{equation*}
En vertu de notre hypothèse, $K(x,y)$ est un polynôme (ou une série bivariée si les sauts ne sont pas bornés). \'Etant donnée une fonction harmonique générique $h=(h(i,j))_{i,j\geq 1}$ avec condition de Dirichlet (i.e., $h(i,0)=0$ pour $i\geq 1$ et $h(0,j)=0$ pour $j\geq 1$), sa \textit{série génératrice} est définie par
\begin{equation*}
   H(x,y)=\sum_{i,j\geq1}h(i,j)x^{i-1}y^{j-1}.
\end{equation*}
L'équation d'harmonicité \eqref{eq:harmonicite} pour $h$ (pour $\lambda=1$) se traduit alors par l'équation fondamentale
\begin{equation}
\label{eq:main_functional_equation}
    K(x,y)H(x,y)=K(x,0)H(x,0) + K(0,y)H(0,y)-K(0,0)H(0,0).
\end{equation}
Montrons à titre d'exemple que l'équation \eqref{eq:main_functional_equation} a bien lieu dans le cas de la marche dite simple: $f(1,0)=f(0,-1)=f(-1,0)=f(0,1)=\frac14$. L'unique fonction harmonique positive vaut $h(i,j) = ij$ et a pour série génératrice $H(x,y)=\frac{1}{(1-x)^2(1-y)^2}$. L'égalité \eqref{eq:main_functional_equation}  provient alors de calculs élémentaires. La série génératrice 
\begin{equation*}
   H(x,y)=\frac{6(xy-1)(y-x)}{(1-x)^4(1-y)^4},
\end{equation*}
qui correspond à la fonction harmonique signée $h(i,j)=ij(i^2-j^2)$, fournit un autre exemple de solution à l'équation \eqref{eq:main_functional_equation}.

\subsection{... et de nombreuses questions}
\label{subsec:ques}

L'équation \eqref{eq:main_functional_equation} montre que toute fonction génératrice d'une fonction harmonique doit prendre la forme (avec des notations évidentes) 
\begin{equation}
\label{eq:FGH}
   H(x,y)=\frac{F(x)+G(y)}{K(x,y)}
\end{equation}
pour des séries entières $F$ et $G$. En revanche, si $K(0,0)=0$, la formule \eqref{eq:FGH} ne définit pas nécessairement une série entière bivariée. Quelles sont les séries $F$ et $G$ telles que \eqref{eq:FGH} est bien une série en $(0,0)$? Quel est le domaine d'analyticité de $H(x,y)$ dans \eqref{eq:FGH}, ou de façon analogue, quelle est la croissance des fonctions harmoniques? Quel choix de $F$ et $G$ conduit à des fonctions harmoniques positives ?

\begin{figure}
\begin{center}
\includegraphics[width=3.5cm]{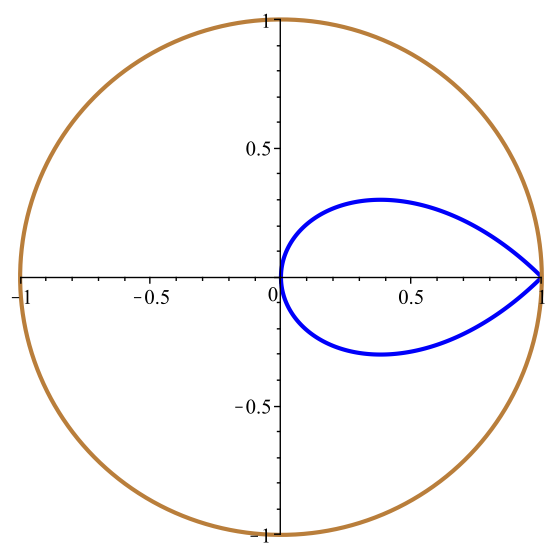}\qquad\qquad
\includegraphics[width=3.5cm]{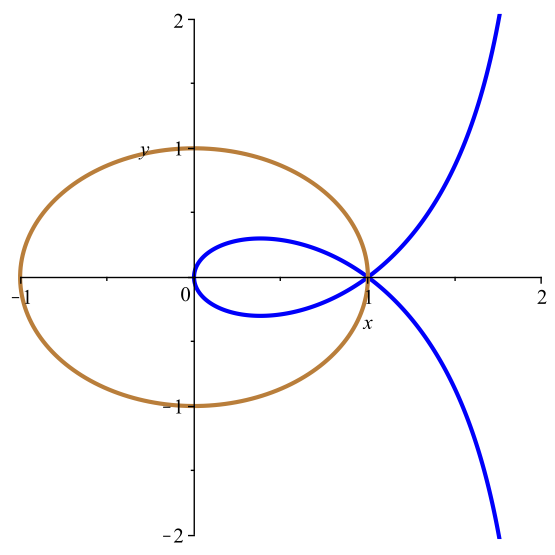}\qquad\qquad
\includegraphics[width=3.5cm]{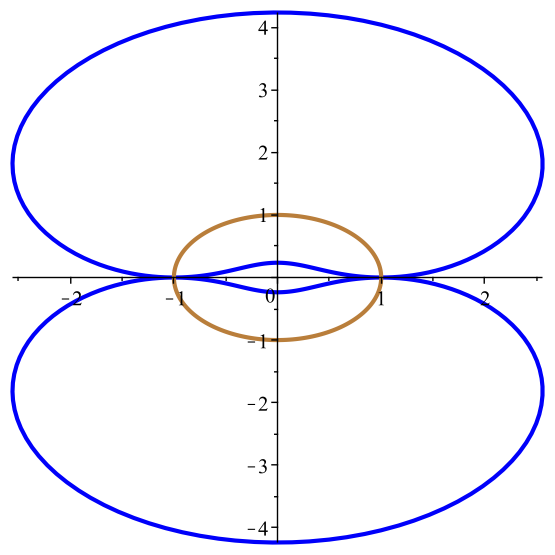}
\end{center}
\caption{Projection en $x$ de l'ensemble \eqref{eq:courbe_Cohen_Boxma}}
\label{fig:some_curves}
\end{figure}

L'ouvrage \cite{CoBo-83} contient une idée féconde pour aborder les questions ci-dessus: dans un contexte probabiliste lié aux distributions stationnaires de marches aléatoires, les auteurs proposent d'évaluer l'équation fonctionnelle principale \eqref{eq:main_functional_equation} sur
\begin{equation}
\label{eq:courbe_Cohen_Boxma}
   \mathcal{K}=\{(x,y)\in\mathbb C^2: \vert x\vert = \vert y \vert\leq 1 \text{ et } K(x,y)=0\}.
\end{equation}
Si, au premier abord peut-être, l'ensemble \eqref{eq:courbe_Cohen_Boxma} paraît surprenant, il présente l'avantage de ne contenir que des $x$ et $y$ petits en module, donc on a de bonnes chances de pouvoir y évaluer nos séries génératrices ; de plus, lors de l'évaluation le terme de gauche de \eqref{eq:main_functional_equation} disparaîtra, comme par définition $K(x,y)=0$. En outre, les rôles joués par $x$ et $y$ sont clairement symétriques dans \eqref{eq:courbe_Cohen_Boxma} ; finalement, les projections de cet ensemble s'avèrent décrire des courbes intéressantes $\mathcal S_x$ et $\mathcal S_y$, voir la Figure~\ref{fig:some_curves}.

À gauche sur la Figure~\ref{fig:some_curves}, en bleu, la courbe $\mathcal S_x=\mathcal S_y$ définie par \eqref{eq:courbe_Cohen_Boxma} dans le cas de la marche simple. Elle est contenue dans la courbe avec auto-intersection représentée au milieu. À droite, la courbe $\mathcal S_x$ pour la marche $(-1,1),(1,0),(0,-1)$ est l'intersection de la courbe représentée avec le disque unité. Le cercle unité est représenté en couleur or (sur les trois graphiques).

\subsection{Réduction à un problème frontière}

Si $h$ est effectivement harmonique, l'évaluation de l'équation \eqref{eq:main_functional_equation} pour $(x,y)\in \mathcal{K}$ donne immédiatement
\begin{equation}
\label{eq:before_BVP}
   F(x)+G(y)=0,
\end{equation}
où $F$ et $G$ sont définies en \eqref{eq:FGH}. Supposons pour la fin de la Section~\ref{sec:analytique} nous trouver dans le cas symétrique $f(k,\ell)=f(\ell,k)$, car alors les projections $\mathcal S_x=\mathcal S_y$ de \eqref{eq:courbe_Cohen_Boxma} sont les mêmes en $x$ et en $y$, ce qui simplifiera notre présentation.

Si $\Psi_{1}$ (resp.\ $\Psi_{2}$) est une \textit{application conforme} entre le demi-plan supérieur $\mathcal H_+=\{z\in\mathbb C: \Im z>0\}$ (resp.\ inférieur $\mathcal H_-$) et le domaine intérieur à $\mathcal S_x=\mathcal S_y$, alors on reformule \eqref{eq:before_BVP} comme le \textit{problème frontière} suivant. Définissons d'abord la fonction 
\begin{equation}
\label{eq:fonction_sec}
   f(t)=\left\{\begin{array}{rcc}
F(\Psi_{1}(t))& \text{si}& t\in\mathcal H_+,\\
-G(\Psi_{2}(t))& \text{si}& t\in\mathcal H_-,
\end{array}\right.
\end{equation}
et faisons tout de suite une première remarque. Supposant $F$ et $G$ analytiques à l'intérieur de $\mathcal S_x$, la fonction $f$ est alors analytique sur $\mathbb C\setminus \mathbb R$. La fonction $f$ est alors dite sectionnellement analytique.

Plus intéressant encore, la condition \eqref{eq:before_BVP} se réécrit simplement comme
\begin{equation}
\label{eq:condition_bord}
   f^+(t)-f^-(t)=0,\quad t\in \mathbb{R}
\end{equation}
où $f^+(t)$ et $f^-(t)$ représentent les limites de $f(z)$ quand $z\to t$ dans $\mathcal H^+$ et $\mathcal H^-$, respectivement. Il est particulièrement frappant qu'après l'utilisation des séries génératrices et différentes reformulations, l'équation d'harmonicité \eqref{eq:harmonicite} devient équivalente à une condition tout à fait classique de continuité \eqref{eq:condition_bord} pour la fonction $f$ en \eqref{eq:fonction_sec}.

Toute fonction entière, en particulier tout polynôme, est continue au voisinage de $\mathbb{R}$ et donc résout \eqref{eq:condition_bord}. Via les changements de variables inverses (notamment dans les applications conformes), on décrit l'ensemble de toutes les fonctions harmoniques. Cette approche permet de répondre aux différentes questions posées dans la Section~\ref{subsec:ques}. Nous conjecturons \cite{HoRaTa-20} que dans le cas d'une dérive nulle, l'unique fonction harmonique positive correspond à prendre $f(t)=t$ dans \eqref{eq:fonction_sec} et \eqref{eq:condition_bord}. Positivité (de la fonction harmonique) rimerait donc avec minimalité (du degré du polynôme $t$ dans $f(t)=t$). 

\section{Applications en combinatoire}

\subsection{Motivations et reformulation probabiliste}

Depuis les années 2000, la communauté combinatoire s'est également emparée de questions étroitement liées aux marches (aléatoires ou déterministes) dans des cônes. Une des motivations principales est l'existence de \textit{bijections} entre des modèles de marches dans des cônes et de nombreux \textit{objets combinatoires}: permutations, tableaux de Young, marches ordonnées, orientations bipolaires, etc., voir par exemple \cite{BMMi-10}. En parallèle de ces motivations exogènes, les marches dans des cônes sont peu à peu devenues un objet combinatoire d'intérêt propre; l'un des objectifs majeurs est de \textit{classifier} les différents modèles de marches selon la complexité de la fonction génératrice de comptage associée, nous y reviendrons plus bas.

Du point de vue combinatoire, le problème typique est le suivant: on fixe un cône $\mathcal C$ de $\mathbb R^d$ et un ensemble de pas ou directions possibles $\mathcal S$ de $\mathbb Z^d$; on s'intéresse alors aux nombres entiers
\begin{equation*}
   c(x;n) \quad \text{et} \quad  c(x,y;n),
\end{equation*}
qui comptent respectivement le nombre de chemins construits à partir de $n$ sauts de $\mathcal S$, partant de $x$ et restant dans le cône pour le premier, et un raffinement de la quantité précédente prenant en compte une arrivée fixée en $y$, de telle sorte que $\sum_{y\in\mathcal C} c(x,y;n)=c(x;n)$.

Le lien avec les quantités probabilistes précédemment introduites est clair: en renormalisant par le cardinal de l'ensemble de pas à la puissance $n$ (qui trivialement compte le nombre de marches partant de $x$ sans contrainte de cône ni de point d'arrivée), on retrouve respectivement la probabilité de survie \eqref{eq:tempslong_discret} et la probabilité locale \eqref{eq:TLL}:
\begin{equation}
\label{eq:interpretation_probab}
   \left\{\begin{array}{rcl}
   \frac{ c(x;n) }{\vert \mathcal S\vert^n} &=& \mathbb{P}(\tau_{x}^{\mathcal C}>n),\medskip\\
   \frac{ c(x,y;n) }{\vert \mathcal S\vert^n} &=& \mathbb P(x+S_n=y,\tau_x^{\mathcal C}>n).
   \end{array}\right.
\end{equation}
\'{E}tant basée sur l'approximation de la marche aléatoire par un mouvement brownien (Section~\ref{sec:derive_nulle}), pour des raisons intrinsèques l'approche probabiliste ne permet pas d'obtenir des expressions exactes pour les nombres combinatoires d'intérêt $c(x;n)$ ou $c(x,y;n)$. Elle est en revanche parfaitement adaptée à des considérations asymptotiques, via les résultats obtenus par Denisov et Wachtel (voir \eqref{eq:tempslong_discret} pour la probabilité de survie et \eqref{eq:TLL} pour la probabilité locale). 

L'approche probabiliste permet aussi de compter (asymptotiquement) des \textit{chemins pondérés} dans des cônes, ce qui par définition revient à poser des probabilités de transition non uniformes sur les pas de $\mathcal S$. C'est une généralisation naturelle; par exemple, certains modèles de marches uniformes en dimension $d\geq 2$ sont équivalents, après projection, à des modèles de marches pondérées en dimension $d'<d$.

En conclusion, la connaissance des fonctions harmoniques discrètes associées aux marches ainsi que de l'exposant critique du cône conduit à des asymptotiques précises pour les nombres de marches $c(x;n)$ et $c(x,y;n)$, via \eqref{eq:tempslong_discret} et \eqref{eq:interpretation_probab}.

\subsection{Classes de fonctions}

De façon plus inattendue, l'interprétation probabiliste \eqref{eq:interpretation_probab} permet d'aborder la question de la \textit{classification} des fonctions génératrices
\begin{equation}
\label{eq:series_generatrices_comb}
   \sum_{n\geq 0} c(x;n)t^n \quad \text{ou} \quad \sum_{n\geq 0} c(x,y;n)t^n
\end{equation}
au sens suivant. \'Etant donnée une série $\sum_{n\geq 0}c_n t^n$, on voudrait la situer dans la hiérarchie suivante des fonctions: 
\begin{itemize}
   \item[$\bullet$]fonctions \textit{rationnelles} (quotients de deux polynômes), 
   \item[$\bullet$]fonctions \textit{algébriques} (solutions d'une équation polynomiale), 
   \item[$\bullet$]fonctions \textit{différentiellement finies} (solutions d'une équation différentielle linéaire à coefficients polynomiaux), 
   \item[$\bullet$]fonctions \textit{différentiellement algébriques} (également solutions d'une équation différentielle à coefficients polynomiaux, mais \textit{a priori} non linéaire).
\end{itemize}
Hors de ces cas avec structure restent les
\begin{itemize}   
   \item[$\bullet$]fonctions \textit{hypertranscendantes} (solutions d'aucune équation différentielle algébrique, à l'instar de la fonction $\Gamma$ d'Euler ou de la fonction $\zeta$ de Riemann).
\end{itemize}
Dire si la série génératrice d'un modèle combinatoire appartient à l'un ou l'autre de ces ensembles de fonctions permet de classifier ce modèle au sein d'une classe plus vaste (par exemple les marches à petits pas dans le quadrant comme dans \cite{BMMi-10,Dretal-18}); cela donne en outre une mesure de sa complexité. Observons également qu'il est équivalent de dire qu'une série est différentiellement finie ou que ses coefficients sont P-récursifs, c'est-à-dire satisfont une équation de récurrence linéaire à coefficients polynomiaux (on peut penser aux nombres de Fibonacci, aux nombres de Catalan ou encore à la factorielle). L'existence de telles récurrences est très naturelle en combinatoire, dans la mesure où elles donnent accès au calcul efficace et rapide de coefficients associés à de grands indices. Ajoutons que les différents ensembles de fonctions susmentionnés satisfont différents théorèmes de clôture (ils restent stables par somme, produit, parfois composition, produit de Hadamard, etc.), ce qui les rend manipulables.

\subsection{Rationalité des exposants critiques et nature des séries génératrices}

Connaître l'asymptotique des coefficients $c_n$ donne une indication, voire une preuve, du fait que leur série génératrice $\sum_{n\geq 0}c_n t^n$ n'appartient pas à l'une des classes de fonctions mentionnées plus haut. Ainsi par exemple, si
\begin{equation}
\label{eq:Ansatz_coeff}
   c_n\sim K \cdot r^n \cdot n^\alpha 
\end{equation}
et si $r$ est un nombre transcendant ou si l'exposant critique $\alpha$ est irrationnel ou un nombre entier négatif, alors la série ne peut pas être algébrique. Dans un esprit similaire, si la série est différentiellement finie, alors nécessairement $r$ et $\alpha$ dans \eqref{eq:Ansatz_coeff} sont tous deux des nombres algébriques.

Dans un contexte plus combinatoire, mentionnons le résultat profond suivant dû à André, Chudnovski et Katz (présenté en détail dans \cite{BoRaSa-14}). Si une suite d'entiers $(c_n)_{n\geq 0}$ se comporte asymptotiquement comme \eqref{eq:Ansatz_coeff} et si la série associée est une nouvelle fois différentiellement finie, alors nécessairement $\alpha$ est rationnel.

Avec cela en tête, il est naturel de se pencher sur la question suivante: étant donné nos nombres de marches $c(x;n)$ et $c(x,y;n)$ et leurs séries associées \eqref{eq:series_generatrices_comb}, sous réserve qu'ils puissent être décrits asymptotiquement comme \eqref{eq:Ansatz_coeff}, quand leur exposant asymptotique est-il rationnel ou irrationnel ? Par exemple, est-il possible d'étudier la rationalité de $p$ dans \eqref{eq:tempslong_discret}?

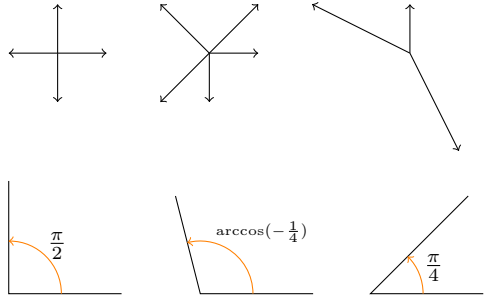
\begin{figure}
\begin{center}
\begin{tikzpicture}[scale=.65] 
    \draw[->,white] (1,2) -- (0,-2);
    \draw[->,white] (1,-2) -- (0,2);
    \draw[->] (0,0) -- (1,0);
    \draw[->] (0,0) -- (-1,0);
    \draw[->] (0,0) -- (0,-1);
    \draw[->] (0,0) -- (0,1);
   \end{tikzpicture}\qquad
\begin{tikzpicture}[scale=.65] 
    \draw[->,white] (1,2) -- (0,-2);
    \draw[->,white] (1,-2) -- (0,2);
    \draw[->] (0,0) -- (1,1);
    \draw[->] (0,0) -- (1,0);
    \draw[->] (0,0) -- (0,-1);
    \draw[->] (0,0) -- (-1,1);
    \draw[->] (0,0) -- (-1,-1);
   \end{tikzpicture}\qquad
\begin{tikzpicture}[scale=.65] 
  \draw[->,white] (1.5,2) -- (-2,-2);
    \draw[->,white] (1.5,-2) -- (-2,2);
        \draw[->] (0,0) -- (0,1);
    \draw[->] (0,0) -- (-2,1);
    \draw[->] (0,0) -- (1,-2);
   \end{tikzpicture}
  \end{center}
    \begin{center}
  \begin{tikzpicture}
  \draw
    (1.5,0) coordinate (a) 
    -- (0,0) coordinate (b) 
    -- (0,1.5) coordinate (c) 
    pic["$\frac{\pi}{2}$", draw=orange, ->, angle eccentricity=1.3, angle radius=0.7cm]
    {angle=a--b--c};
\end{tikzpicture}\qquad
\begin{tikzpicture}
  \draw
    (1.5,0) coordinate (a) 
    -- (0,0) coordinate (b) 
    -- (-0.33,1.30) coordinate (c) 
    pic["\quad\tiny{$\arccos\bigl(-\frac{1}{4}\bigr)$}", draw=orange, ->, angle eccentricity=1.5, angle radius=0.7cm]
    {angle=a--b--c};
\end{tikzpicture}\qquad
\begin{tikzpicture}
  \draw
    (1.5,0) coordinate (a) 
    -- (0,0) coordinate (b) 
    -- (1.3,1.3) coordinate (c) 
    pic["$\frac{\pi}{4}$", draw=orange, ->, angle eccentricity=1.3, angle radius=0.7cm]
    {angle=a--b--c};
\end{tikzpicture}
\end{center}
\caption{Exemples en dimension 2: des modèles de marches et la valeur de $\arccos a$ dans \eqref{eq:arccos} associée}
\label{fig:dim-2}
\end{figure}

Cela ouvre tout un champ de questionnements. Il y a peu à dire en dimension $1$, car les exposants typiques sont rationnels (le plus souvent $0$, $-\frac{1}{2}$, $-\frac{3}{2}$). En dimension deux, on utilise d'abord \eqref{eq:tempslong_discret} et \eqref{eq:interpretation_probab} pour voir que l'exposant critique $\alpha$ dans \eqref{eq:Ansatz_coeff} vaut $\alpha=-\frac{p}{2}$ pour $c(x;n)$ en dérive nulle. Un théorème de la limite locale (asymptotique de \eqref{eq:TLL}) et \eqref{eq:interpretation_probab} prouveraient que $\alpha=-p-1$ pour $c(x,y;n)$ sans condition sur la dérive. Les résultats de Denisov et Wachtel conduisent à
\begin{equation}
\label{eq:arccos}
   p=\frac{\pi}{\arccos a},
\end{equation}
où $a$ est naturellement associé à la marche aléatoire et s'interprète comme un coefficient de corrélation. L'algébricité de $a$ est immédiate par construction (il est calculé en résolvant un système d'équations polynomiales). Mieux, la rationalité de $p$ dans \eqref{eq:arccos} est accessible et permet (voir \cite{BoRaSa-14}) de prouver que tout un ensemble de modèles admet des séries non différentiellement finies. Ainsi les premier et dernier modèles sur la Figure~\ref{fig:dim-2} possèdent un $\alpha$ rationnel, tandis qu'on prouve facilement que $p=\frac{\pi}{\arccos\left( -\frac{1}{4}\right)}\notin \mathbb Q$ pour la marche du milieu sur la Figure~\ref{fig:dim-2}. 

\begin{figure}
\begin{center}
\includegraphics[width=3.5cm]{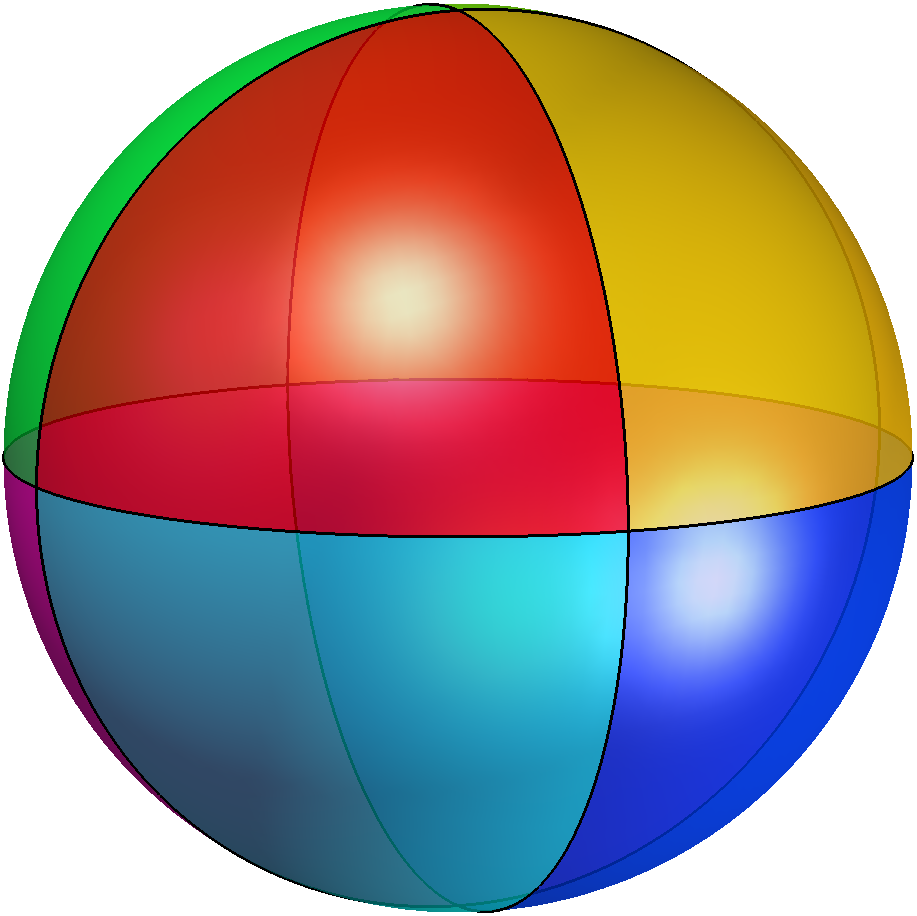}\qquad\qquad\qquad
\includegraphics[width=3.5cm]{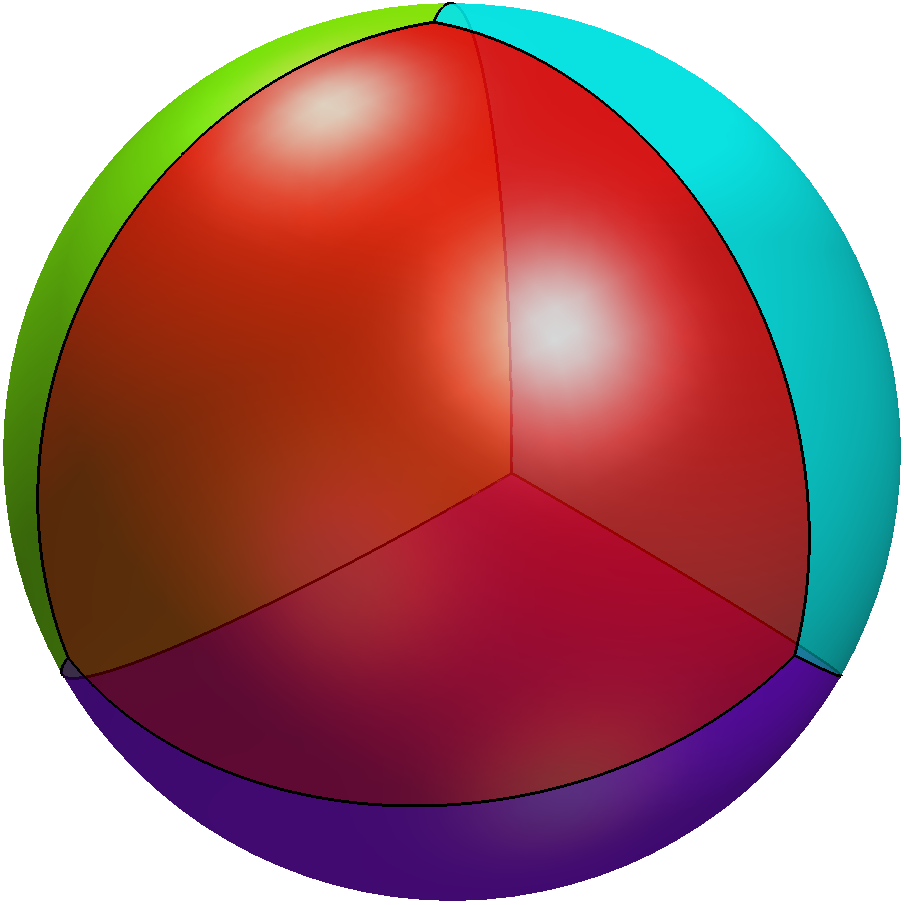}
\end{center}
\caption{En dimension $3$ l'exposant critique est li\'e \`a la valeur propre principale pour le probl\`eme de Dirichlet sur un triangle sph\'erique}
\label{fig:dim-3}
\end{figure}

En dimension $d\geq 3$, Denisov et Wachtel \cite{DeWa-15} expriment $\alpha$ au moyen de la valeur propre principale d'un problème de Dirichlet. De ce fait, on se heurte dès la dimension $3$ à des questions de théorie spectrale encore ouvertes, comme par exemple le problème suivant : on considère un triangle sphérique et on essaie de calculer la première valeur propre pour le problème de Dirichlet. \`A gauche sur la Figure~\ref{fig:dim-3}, la marche simple dans l'octant correspond à un triangle sphérique équilatéral avec angle $\frac{\pi}{2}$, dont le spectre est parfaitement connu. \`A droite sur la même figure, un grand triangle équilatéral d'angle $\frac{2\pi}{3}$, qui correspond à une marche dans l'octant $\mathbb N^3$ avec sauts égaux aux trois vecteurs de la base canonique plus $(-1,-1,-1)$. Les combinatoriciens nomment cet ensemble de pas "Kreweras 3D", en lien avec les anciens travaux de Germain Kreweras. On conjecture que le premier exposant est irrationnel dans ce second cas.

\subsection*{Remerciements}
Nous remercions chaleureusement Irina Ignatiouk-Robert et Sami Mustapha pour leur disponibilité et leurs remarques. Nous adressons également nos vifs remerciements à Mireille Bousquet-Mélou, Aurélien Djament et Mylène Maïda pour leur relecture attentive de l'article et leurs nombreuses suggestions d'amélioration. 

\bibliographystyle{abbrv}

\medskip

\noindent Kilian Raschel \href{mailto:raschel@math.cnrs.fr}{\texttt{raschel@math.cnrs.fr}} 

\noindent Universit\'e d'Angers, CNRS, Laboratoire Angevin de Recherche en Math\'ematiques, SFR MATHSTIC, Angers, France

\medskip

\noindent Pierre Tarrago \href{mailto:pierre.tarrago@sorbonne-universite.fr}{\texttt{pierre.tarrago@sorbonne-universite.fr}} 

\noindent Sorbonne Université, Laboratoire de Probabilités, Statistique et Modélisation, Paris, France 

\medskip

\noindent This project has received funding from the European Research Council (ERC) under the European Union's Horizon 2020 research and innovation programme under the Grant Agreement No.\ 759702 and from Centre Henri Lebesgue, programme ANR-11-LABX-0020-0

\end{document}